# STABILITY OF SOLUTIONS OF FUZZY DIFFERENTIAL EQUATIONS

LE VAN HIEN

**Abstract.** In this paper, we study the stability of solutions of fuzzy differential equations by Lyapunov's second method. By using scale equations and comparison principle for Lyapunov - like functions, we give some sufficient criterias for the stability and asymptotic stability of solutions of fuzzy differential equations.

## 1 Introductions

Fuzzy sets were introduced in 1965 by Lotfi Zadeh [20] with a view to reconsile mathematical modeling and human knowledge in the engineering science. Since then, a considerable body of literature has blossomed around the concepts of fuzzy sets in an incredibly wide range of areas, from mathematics and logics to traditional and advanced engineering methodologies (from civil engineering to computational intelligence). Applications are found in many contexts, from medicine to finance, from human factors to consumer products...Fuzzy logic is now currently used in the industrial practice of advanced information technology.

Recently, the industrial interest in fuzzy control and logic [3] has dramatically increased the study of fuzzy systems. The calculus of fuzzy - valued functions has been initiated [4, 5, 6, 9, 16] and the study of the initial value problems for fuzzy differential equations has been initiated in [1, 2, 7, 8, 12, 13]. The existence and uniqueness of solutions of fuzzy differential equations is considered by some authors [12, 13, 14, 18]. The extension of solutions were given in [21] and the global existence of solutions were given in [17].

The investigation of stability of solutions is the most important problem in the qualitative theory of differential equations. It has been widely applied in Physic, Mechanic, Control,...

In this paper, we study the stability theory which corresponds to Lyapunov stability theory for fuzzy differential equations. By using differential inequalities and comparison principle, somes sufficient conditions for the stability of solutions of fuzzy differential equations were estimated.

## 2 Preliminaries

Let $P_K(\mathbb{R}^n)$ denotes the familiy of all non-empty compact, convex subsets of $\mathbb{R}^n$ and define the addition and scalar multiplication in $P_K(\mathbb{R}^n)$ as usual. Let $A, B$ are two non-empty subsets in $\mathbb{R}^n$.

The distance between $A$ and $B$ is defined by Haussdorff metric:

$$d_H(A, B) = \max[\sup_{a \in A} \inf_{b \in B} \|a - b\|, \sup_{b \in B} \inf_{a \in A} \|a - b\|]$$

where $\|.\|$ denotes a norm in $\mathbb{R}^n$. Then it is clear that $(P_K(\mathbb{R}^n), d_H)$ becomes a metric space. Moreover, the metric space $(P_K(\mathbb{R}^n), d_H)$ is complete and separable (see [15]). Let $T = [a; b], a \geq 0$ be a interval in $\mathbb{R}$ and denote $\varepsilon^n = \{u : \mathbb{R}^n \longrightarrow [0; 1] | u$ satisfies (i) to (iv) below $\}$.

(i) $u$ is normal, that is, there exists $x_0 \in \mathbb{R}^n$ such that $u(x_0) = 1$;

(ii) $u$ is fuzzy convex, that is, for $x, y \in \mathbb{R}^n$ and $0 \leq \lambda \leq 1$:

$$u(\lambda x + (1 - \lambda) y) \geq \min[u(x), u(y)];$$

(iii) $u$ is upper semi-continuous;

(iv) $[u]^0 = \overline{\{x \in \mathbb{R}^n : u(x) > 0\}}$ is compact subset in $\mathbb{R}^n$;

For $0 < \alpha \leq 1$, we denote $[u]^\alpha = \{x \in \mathbb{R}^n : u(x) \geq \alpha\}$, then from $(i)$ to $(iv)$, it follow that the $\alpha$−level $[u]^\alpha \in P_K(\mathbb{R}^n)$ for all $\alpha \in [0; 1]$. For later purpose, we define $\widehat{o} \in \varepsilon^n$ as $\widehat{o}(x) = \chi_{\{0\}}(x) = 1$ if $x = 0$ and $\widehat{o}(x) = 0$ if $x \neq 0$. Define a metric function $d : \varepsilon^n \times \varepsilon^n \longrightarrow \mathbb{R}_+$ by

$$d[u, v] = \sup_{0 \leq \alpha \leq 1} d_H([u]^\alpha, [v]^\alpha)$$

then $(\varepsilon^n, d)$ becomes a complete metric space (see [12, 15]). We list here some properties of metric $d[u, v]$ (see [7, 12, 14, 15]).

(i) $d[u, w] \leq d[u, v] + d[v, w]$;

(ii) $d[\lambda u, \lambda v] = |\lambda| d[u, v]$;

(iii) $d[u + w, v + w] = d[u, v], u, v, w \in \varepsilon^n, \lambda \in \mathbb{R}$;

For $x, y \in \varepsilon^n$, if there exists $z \in \varepsilon^n$ such that $x = y + z$ then $z$ is called H-difference of $x$ and $y$ and is denoted by $x - y$.

A mapping $F : T \longrightarrow \varepsilon^n$ is differentiable at $t_0 \in T$ if for small $h > 0$, there exist H-differences $F(t_0 + h) - F(t_0); F(t_0) - F(t_0 - h)$ and there exists a $F'(t_0) \in \varepsilon^n$ such that the limits

$$\lim_{h \to 0^+} \frac{F(t_0 + h) - F(t_0)}{h}; \lim_{h \to 0^+} \frac{F(t_0) - F(t_0 - h)}{h}$$

exist and equal $F'(t_0)$. If $F, G$ differentiable at $t$ then $(F + G)'(t) = F'(t) + G'(t)$ and $(\lambda F)'(t) = \lambda F'(t), \lambda \in \mathbb{R}$ (see [7, 12, 16]).

If $F : T \longrightarrow \varepsilon^n$ is strongly measurable and integrably bounded then it is integrable on $T$ and $\int_T F(t) dt \in \varepsilon^n$,

$$[\int_T F(t) dt]^\alpha = \int_T F_\alpha(t) dt, 0 < \alpha \leq 1; F_\alpha(t) = [F(t)]^\alpha.$$

where $\int_T F_\alpha(t) dt$ is *Aumann integral*. It is wellknown that $[\int_T F(t) dt]^0 = \int_T F_0(t) dt$ (see [7], Remark 4.1). Also the following properties of integral are valid (see [5, 6, 7, 12]). If $F, G : T \longrightarrow \varepsilon^n$ be integrable on $T$ and $\lambda \in \mathbb{R}$ then:

(i) $\int_T (F+G)(t)dt = \int_T F(t)dt + \int_T G(t)dt$;

(ii) $\int_T (\lambda F)(t)dt = \lambda \int_T F(t)dt$;

(iii) $d[F(.), G(.)] : T \longrightarrow \mathbb{R}_+$ is integrable;

(iv) $d[\int_T F(t)dt, \int_T G(t)dt] \leq \int_T d[F(t), G(t)]dt$;

(v) $\int_a^b F(t)dt = \int_a^c F(t)dt + \int_c^b F(t)dt, a \leq c \leq b$;

If $F$ is continuous then $G(t) = \int_a^t F(\tau)d\tau$ is differentiable on $T$ and $G'(t) = F(t), \forall t \in T$. Moreover, if $F$ is differentiable on $T$ and $F'(.)$ is integrable on $T$ then for all $t \in T$ we have $F(t) = F(t_0) + \int_{t_0}^t F(\tau)d\tau, a \leq t_0 \leq t \leq b$. If $F$ is continuous on $T$ and $G(t) = \int_a^t F(\tau)d\tau$ then for $t_1 \leq t_2$ we have (see [7])

$$d[G(t_1), G(t_2)] \leq (t_2 - t_1) \sup_{[t_1, t_2]} d[F(t), \widehat{o}]$$

## 3 Stability

Consider fuzzy differential equation:

$$\frac{dx}{dt} = f(t, x), x(t_0) = x_0 \tag{3.1}$$

where $f \in C[\mathbb{R}_+ \times S(\rho), \varepsilon^n], S(\rho) = \{x \in \varepsilon^n : d[x, \widehat{o}] < \rho\}, f(t, \widehat{o}) \equiv \widehat{o}$. Hence, equation (3.1) has trivial solution $x = \widehat{o}$.

In this section, we shall discuss the stability, especially, asymptotically stability of solutions of Eq(3.1) by Lyapunov's second method. First, we give some notions of stability which are used in the sequel. Let $x(t) = x(t; t_0, x_0)$ be any solution of (3.1) existing on $[t_0, \infty)$. Denote $\mathcal{K} = \{a \in C[\mathbb{R}_+, \mathbb{R}_+], a(0) = 0, a(.) \text{ is increasing}\}$.

**Definition 1.** The trivial solution $x = \widehat{o}$ of (3.1) is stable if for any $\varepsilon > 0, t_0 \in \mathbb{R}_+$, there exists a $\delta = \delta(t_0, \varepsilon) > 0$ such that if $d[x_0, \widehat{o}] < \delta$ then $d[x(t), \widehat{o}] < \varepsilon, \forall t \geq t_0$.

If the $\delta$ in the above definition is independent of $t_0$ then $x = \widehat{o}$ is said to be uniformly stable.

**Definition 2.** The trivial solution $x = \widehat{o}$ of (3.1) is asymptotically stable if $x = \widehat{o}$ is stable and for any $t_0 \in \mathbb{R}_+$, there exists a $\Delta = \Delta(t_0) > 0$ such that if $d[x_0, \widehat{o}] < \Delta$ then $\lim_{t \to \infty} d[x(t; t_0, x_0), \widehat{o}] = 0$.

**Definition 3.** The trivial solution $x = \widehat{o}$ of (3.1) is uniformly asymptotically stable if it is uniformly stable and there exists a $\delta_0 > 0$ such that for any $\varepsilon > 0$, there exist $T(\varepsilon) \geq 0$ such that if $d[x_0, \widehat{o}] < \delta_0, t_0 \in \mathbb{R}_+$ then

$$d[x(t; t_0, x_0), \widehat{o}] < \varepsilon, \forall t \geq t_0 + T(\varepsilon)$$

**Definition 4.** The trivial solution $x = \widehat{o}$ is exponentially stable if any solution $x(t) = x(t; t_0, x_0)$ of (3.1) satisfies:

$$d[x(t), \widehat{o}] \leq \beta(d[x_0, \widehat{o}], t_0)e^{-\alpha(t - t_0)}, t \geq t_0$$

where $\beta(h, t) : [0, H) \times \mathbb{R}_+ \longrightarrow \mathbb{R}_+$ increasing in $h \in [0, H)$ for some $H > 0$ and $\alpha$ is a positive constant. If $H = \infty$ then $x = \widehat{o}$ is called global exponentially stable.

If the function $\beta(., .)$ does not depend on $t_0$ then $x = \widehat{o}$ is called uniformly exponentially stable.

Before prove the stability of solutions of (3.1), we need the following Lemma (see [11, 19] for details).

**Lemma 3.1.** *Let $g(t,x)$ be a continuous function on $\mathbb{R}_+^2$ and $r(t) = r(t; t_0, w_0), r(t_0) = w_0$ be the maximal solution of the scalar differential equation:*

$$w' = g(t, w) \tag{3.2}$$

*existing on $[t_0, \infty)$. Let $m(t)$ be a continuous function on $\mathbb{R}_+$ satisfies*

$$d^+ m(t) = \limsup_{h \to 0^+} \frac{m(t+h) - m(t)}{h} \leq g(t, m(t)), t \geq t_0$$

*Then $m(t) \leq r(t), \forall t \geq t_0$ if $m(t_0) \leq w_0$.*

Let $V(t, x) : \mathbb{R}_+ \times S(\rho) \longrightarrow \mathbb{R}$ be a given function. Then we define

$$D_f^+ V(t, x) = \limsup_{h \to 0^+} \frac{1}{h}[V(t+h, x + hf(t, x)) - V(t, x)]$$

where $f(.)$ is the right-hand side of (3.1) and $D_f^+ V$ is called the upper derivation of $V(t, x)$ along the trajectory of (3.1). Let $x(t)$ be a solution of (3.1) then $d^+ V(t, x(t))$ denotes the upper derivation of $V(t, x(t))$, i.e.

$$d^+ V(t, x(t)) = \limsup_{h \to o^+} \frac{1}{h}[V(t+h, x(t+h)) - V(t, x(t))]$$

Note that, if $V(t, x)$ is Lipchitzian in $x$ then we have $d^+ V(t, x(t)) \leq D_f^+ V(t, x(t))$.

**Theorem 3.1.** *Suppose that there exists a function $V(t, x)$ satisfies the following conditions:*

*(i) $|V(t, x) - V(t, y)| \leq L(t) d[x, y], \forall (t, x), (t, y) \in \mathbb{R}_+ \times S(\rho), L(.) \in C[\mathbb{R}_+, \mathbb{R}_+]$;*

*(ii) $a(d[x, \widehat{o}]) \leq V(t, x), V(t, \widehat{o}) = 0$, where $a(.) \in \mathcal{K}$ class;*

*(iii) $D_f^+ V(t, x) \leq g(t, V(t, x)), g(., .) \in C[\mathbb{R}_+^2, \mathbb{R}]; g(t, 0) = 0$;*

*If the solution $w = 0$ of the equation in the form (3.2) is stable (asymptotically stable) then the trivial solution $x = \widehat{o}$ of (3.1) is stable (asymptotically stable).*

*Proof.* Let $x(t) = x(t; t_0, x_0), t_0 \in \mathbb{R}_+$ be any solution of Eq(3.1) existing on $[t_0, \infty)$ and solution $w = 0$ of (3.2) is stable. Then, for any $0 < \varepsilon < \rho$, exists a $\delta_0 = \delta_0(t_0, \varepsilon) > 0$ such that if $0 \leq w_0 < \delta_0$ then $|w(t; t_0, w_0)| < a(\varepsilon), \forall t \geq t_0$. From (ii), it follows that there exists $\delta = \delta(t_0, \varepsilon) > 0$ such that $V(t_0, x) < \delta_0$ if $d[x, \widehat{o}] < \delta$. We will show that if $d[x_0, \widehat{o}] < \delta$ then $d[x(t), \widehat{o}] < \varepsilon, \forall t \geq t_0$.

Suppose that $d[x(t), \widehat{o}] \geq \varepsilon$ for some $t_* > t_0$ then there exists a $t_1 > t_0$ such that

$$d[x(t_1), \widehat{o}] = \varepsilon; d[x(t), \widehat{o}] < \varepsilon, \forall t \in [t_0, t_1)$$

Let $m(t) = V(t, x(t)), t \geq t_0$ then we have:

$$\begin{aligned} m(t+h) - m(t) &= V(t+h, x(t+h)) - V(t, x(t)) \\ &= V(t+h, x(t+h)) - V(t+h, x(t) + hf(t, x(t))) + \\ &\quad + V(t+h, x(t) + hf(t, x(t))) - V(t, x(t)) \\ &\leq L(t+h) d[x(t+h), x(t) + hf(t, x(t))] + \\ &\quad + V(t+h, x(t) + hf(t, x(t))) - V(t, x(t)). \end{aligned}$$

For small $h > 0$, H-differences of $x(t+h)$ and $x(t)$ is assumed to exsit. Let $x(t+h) = x(t) + z(t)$ and using the properties of metric $d[x,y]$ we have:

$$d[x(t+h), x(t) + hf(t, x(t))] = hd[\frac{x(t+h) - x(t)}{h}, f(t, x(t))]$$

Hence

$$\begin{aligned}
d^+ m(t) &= \limsup_{h \to 0^+} \frac{1}{h}[m(t+h) - m(t)] \\
&\leq L(t) \limsup_{h \to 0^+} d[\frac{x(t+h) - x(t)}{h}, f(t, x(t))] + \\
&\quad + \limsup_{h \to 0^+} \frac{1}{h}[V(t+h, x(t) + hf(t, x(t))) - V(t, x(t))] = \\
&= L(t)d[x'(t), f(t, x(t))] + D_f^+ V(t, x(t)) = D_f^+ V(t, x(t)) \\
&\leq g(t, m(t)), t_0 \leq t \leq t_1.
\end{aligned}$$

Applying Lemma 3.1, $m(t) \leq r(t; t_0, w_0), w_0 = V(t_0, x_0), t \in [t_0, t_1]$.
On the other hand, $V(t_0, x_0) < \delta_0$, so, $r(t; t_0, w_0) < a(\varepsilon), t \in [t_0, t_1]$ and therefore

$$m(t_1) \leq r(t_1; t_0, w_0) < a(\varepsilon)$$

By the choice of $t_1$, we have $a(\varepsilon) = a(d[x(t_1), \hat{o}]) \leq V(t_1, x(t_1)) = m(t_1) < a(\varepsilon)$. This is a contradiction, hence

$$d[x(t), \hat{o}] < \varepsilon, \forall t \geq t_0$$

This shows that the trivial solution $x = \hat{o}$ of (3.1) is stable.

　　　If $w = 0$ of (3.2) is asymptotically stable then it's stable, therefore $x = \hat{o}$ of (3.1) is stable. For $t_0 \in \mathbb{R}_+$, there exist $\delta = \delta(t_0) > 0, \Delta_1(t_0) > 0$ such that $d[x(t), \hat{o}] < \rho, \forall t \geq t_0$ if $d[x_0, \hat{o}] < \delta$ and if $0 \leq w_0 < \Delta_1(t_0)$ then $\lim_{t \to \infty} w(t; t_0, w_0) = 0$. From hypothesises of function $V(t, x)$, we can find a $\Delta_2 > 0$ such that if $d[x, \hat{o}] < \Delta_2$ then $V(t, x) < \Delta_1(t_0)$. Put $\Delta = \min[\delta, \Delta_2]$. Let $x(t)$ be any solution of (3.1), $t_0 \in \mathbb{R}_+, d[x_0, \hat{o}] < \Delta$. Define $m(t) = V(t, x(t)), t \geq t_0$. By the first part of this proof we see that $d^+ m(t) \leq g(t, m(t))$. Apply Lemma 3.1,

$$m(t) \leq r(t; t_0, w_0), w_0 = V(t_0, x_0), t \geq t_0$$

Since $w_0 = V(t_0, x_0) < \Delta_1(t_0)$, so $\lim_{t \to \infty} r(t; t_0, w_0) = 0$.

From $a(d[x(t), \hat{o}]) \leq V(t, x(t)) = m(t) \leq r(t; t_0, w_0), a(.) \in \mathcal{K}$, it follows that $\lim_{t \to \infty} d[x(t), \hat{o}] = 0$. This shows that $x = \hat{o}$ is asymptotically stable. The proof is completed. □

**Theorem 3.2.** *Suppose that there exists a function $V(t, x)$ satisfies:*

(i) $|V(t, x) - V(t, y)| \leq L(t)d[x, y], \forall (t, x), (t, y) \in \mathbb{R}_+ \times S(\rho), L(.) \in C[\mathbb{R}_+, \mathbb{R}_+]$;

(ii) $a(d[x, \hat{o}]) \leq V(t, x) \leq b(d[x, \hat{o}]), a(.), b(.) \in \mathcal{K}$;

(iii) $D_f^+ V(t, x) \leq g(t, V(t, x)), g \in C[\mathbb{R}_+^2, \mathbb{R}]; g(t, 0) = 0$;

*If the solution $w = 0$ of (3.2) is uniformly stable (uniformly asymptotically stable) then the trivial solution $x = \hat{o}$ of (3.1) is uniformly stable (uniformly asymptotically stable).*

*Proof.* If solution $w = 0$ of (3.2) is uniformly stable then for any $\varepsilon > 0$ there exists $\delta_0 > 0$ such that if $t_0 \in \mathbb{R}_+$ and $0 \leq w_0 < \delta_0$ then $|w(t; t_0, w_0)| < a(\varepsilon), \forall t \geq t_0$. By choosing a $\delta = \delta(\varepsilon) > 0$ such that $b(\delta) < a(\delta_0)$ and by the same argument in the proof of Theorem 3.1, it can be proved that if $d[x_0, \widehat{o}] < \delta$ then $d[x(t; t_0, x_0), \widehat{o}] < \varepsilon, t \geq t_0$. This shows that $x = \widehat{o}$ is uniformly stable.

Now, we assume $w = 0$ is uniformly asymptotically stable, then by the first part of this proof, the trivial solution $x = \widehat{o}$ is uniformly stable. Hence, there exists $\delta_0 > 0$ such that, $t_0 \in \mathbb{R}_+, d[x_0, \widehat{o}] < \delta_0$ implies $d[x(t; t_0, x_0), \widehat{o}] < \rho, \forall t \geq t_0$. Moreover, there exists $\delta_1 > 0$ such that for any $\varepsilon > 0$, exists $T = T(\varepsilon) \geq 0$ such that if $t_0 \geq 0, 0 \leq w_0 < \delta_1$ then $|w(t; t_0, w_0)| < a(\varepsilon), \forall t \geq t_0 + T$. Put $\delta = \min[\delta_0, b^{-1}(\delta_1)]$. By the same argument in the proof of Theorem 3.1, it can be proved that, if $d[x_0, \widehat{o}] < \delta$ then $d[x(t; t_0, x_0), \widehat{o}] < \varepsilon, \forall t \geq t_0 + T(\varepsilon)$. This shows that $x = \widehat{o}$ is uniformly asymptotically stable. The proof is completed. □

**Example 3.1.** *Consider a fuzzy-valued function $f(t, x)$ which satisfies*

$$d[f(t, x), \widehat{o}] \leq a(t) d[x, \widehat{o}]; \int_0^\infty a(t) dt < \infty$$

*(for example $f(t, x) = \dfrac{1}{1 + t^2} x, a(t) = \dfrac{1}{1 + t^2}$ satisfies all the above conditions).*

*Then the trivial solution $x = \widehat{o}$ of (3.1) is uniformly stable.*

*Proof.* Consider a Lyapunov function $V(t, x) = d[x, \widehat{o}]$.
Then $\dfrac{1}{2} d[x, \widehat{o}] \leq V(t, x) \leq 2 d[x, \widehat{o}]$ and $|V(t, x) - V(t, y)| \leq d[x, y], \forall (t, x); (t, y) \in \mathbb{R}_+ \times \varepsilon^n$. For $h > 0$, we have:

$$\begin{aligned} V(t + h, x + h f(t, x)) &= d[x + h f(t, x), \widehat{o}] \\ &\leq d[x, \widehat{o}] + h d[f(t, x), \widehat{o}] \\ &\leq d[x, \widehat{o}] + h a(t) d[x, \widehat{o}] \end{aligned}$$

Hence, $D_f^+ V(t, x) \leq a(t) d[x, \widehat{o}] = g(t, V(t, x))$, where $g(t, w) = a(t) w$. It's easy to show that the solution $w = 0$ of (3.2) is uniformly stable, so by Theorem 3.2, the trivial solution $x = \widehat{o}$ of (3.1) is uniformly stable. □

**Theorem 3.3.** *Suppose that:*

(i) $f(t, x)$ *is bounded on* $\mathbb{R}_+ \times S(\rho)$;

(ii) $\exists V(t, x)$ *satisfies* $|V(t, x) - V(t, y)| \leq L(t) d[x, y]; a(d[x, \widehat{o}]) \leq V(t, x) \leq a_0(t, d[x, \widehat{o}])$, *where* $a(.) \in \mathcal{K}, a_0(t, .) \in \mathcal{K}$ *for each* $t \in \mathbb{R}_+$;

(iii) $D_f^+ V(t, x) + V^*(t, x) \leq g(t, V(t, x))$, *where* $g \in C[\mathbb{R}_+ \times \mathbb{R}, \mathbb{R}], g(t, .)$ *is non-decreasing for each* $t \in \mathbb{R}_+$ *and* $V^* \in C[\mathbb{R}_+ \times S(\rho), \mathbb{R}_+], V^*(t, x) \geq c(d[x, \widehat{o}]), c(.) \in \mathcal{K}$.

*If solution $w = 0$ of (3.2) is stable then the trivial solution $x = \widehat{o}$ of (3.1) is asymptotically stable.*

*Proof.* By Theorem 3.1, the trivial solution $x = \widehat{o}$ is stable. Hence, for $t_0 \in \mathbb{R}_+$, there exists $\delta_1(t_0)$ such that $d[x_0, \widehat{o}] < \delta_1$ implies $d[x(t; t_0, x_0), \widehat{o}] < \rho, \forall t \geq t_0$. Moreover, for $t_0 \in \mathbb{R}_+$, there exists $\delta_2(t_0) > 0$ such that if $0 \leq w_0 < \delta_2(t_0)$ then $|r(t; t_0, w_0)| < \rho, \forall t \geq t_0$, where $r(t; t_0, w_0)$ is the maximal solution of (3.2). Since $a_0(t_0, .) \in \mathcal{K}$, there exists $\delta_3(t_0) > 0$ such that $a_0(t_0, \delta_3) < \delta_2(t_0)$. Put $\delta = \delta(t_0) = \min\{\delta_1, \delta_2, \delta_3\}$. Let $x(t) = x(t; t_0, x_0)$ be any solution of (3.1), $d[x_0, \widehat{o}] < \delta$. We will show that

$$\lim_{t \to \infty} d[x(t), \widehat{o}] = 0$$

Suppose that $\limsup_{t \to \infty} d[x(t), \widehat{o}] > 0$. Then there exists $\eta > 0$ and a sequence $\{t_n\} \to \infty$ such that

$$d[x(t_n), \widehat{o}] \geq \eta, n = 0, 1, 2, \ldots$$

By the boundedness of $f(t, x)$ and by taking a subsequence of $\{t_n\}$, we can assume that there exist $M > 0, \{t_n\} \to \infty$ such that $t_{n+1} - t_n \geq \dfrac{\eta}{2M}, n \geq 0$.

For $t \in [t_n, t_n + \dfrac{\eta}{2M}]$, we have $x(t) = x(t_n) + \int_{t_n}^t f(\tau, x(\tau))d\tau$. Hence

$$d[x(t), \widehat{o}] \geq d[x(t_n), \widehat{o}] - \int_{t_n}^t d[f(\tau, x(\tau)), \widehat{o}]d\tau \geq \eta - M\frac{\eta}{2M} = \frac{\eta}{2}.$$

Define $m(t) = V(t, x(t)) + \int_{t_0}^t V^*(\tau, x(\tau))d\tau, t \geq t_0$. Then:

$$d^+ m(t) \leq D_f^+ V(t, x(t)) + V^*(t, x(t)) \leq g(t, V(t, x(t))) \leq g(t, m(t)), t \geq t_0$$

Applying Lemma 3.1, it follows that $m(t) \leq r(t; t_0, w_0)$, where $w_0 = V(t_0, x_0)$.

Since $V(t_0, x_0) \leq a_0(t_0, d[x_0, \widehat{o}]) < a_0(t_0, \delta) < \delta_2(t_0)$, hence $|r(t; t_0, w_0)| < \rho, \forall t \geq t_0$. Therefore,

$$V(t_n + \frac{\eta}{2M}, x(t_n + \frac{\eta}{2M})) \leq r(t; t_0, w_0) - \sum_{k=0}^n \int_{t_k}^{t_k + \frac{\eta}{2M}} V^*(\tau, x(\tau))d\tau$$

$$\leq r(t; t_0, w_0) - c(\frac{\eta}{2})\frac{\eta}{2M}n$$

$$< \rho - c(\frac{\eta}{2})\frac{\eta}{2M}n < 0$$

for $n$ sufficiently large. This is a contradiction and therefore:

$$\lim_{t \to \infty} d[x(t), \widehat{o}] = 0$$

The proof is completed. □

**Theorem 3.4.** *Let the assumptions (i), (ii) of Theorem 3.2 hold and*

(iii') $D_f^+ V(t, x) + V^*(t, x) \leq g(t, V(t, x)), g(.,.)$ *as in Theorem 3.3,*
$V^* \in C[\mathbb{R}_+ \times S(\rho), \mathbb{R}_+], V^*(t, x) \geq c(d[x, \widehat{o}]), c(.) \in \mathcal{K}.$

*If solution $w = 0$ of (3.2) is uniformly stable then the trivial solution $x = \widehat{o}$ of (3.1) is uniformly asymptotically stable.*

*Proof.* By Theorem 3.2, $x = \widehat{o}$ of (3.1) is uniformly stable. Hence, for $\varepsilon = \rho$ there exists $\delta_0 > 0$ such that if $t_0 \in \mathbb{R}_+, d[x_0, \widehat{o}] < \delta_0$ then

$$d[x(t; t_0, x_0), \widehat{o}] < \rho, t \geq t_0$$

We can assume that $\delta_0$ satisfies if $0 \leq w_0 < b(\delta_0)$ then $|r(t; t_0, w_0)| < a(\rho), \forall t \geq t_0$. By the uniformly stability of $x = \widehat{o}$, for any $\varepsilon > 0$, there exists $\delta > 0$ such that if $d[x_0, \widehat{o}] < \delta, t_0 \in \mathbb{R}_+$ then $d[x(t; t_0, x_0), \widehat{o}] < \varepsilon, \forall t \geq t_0$. Let's putting $T = T(\varepsilon) = 1 + \dfrac{a(\rho)}{c(\delta)}$. Let $x(t) = x(t; t_0, x_0)$ be any solution of (3.1), $d[x_0, \widehat{o}] < \delta_0$. We will show that $d[x(t), \widehat{o}] < \delta$ for some $t_* \in [t_0, t_0 + T(\varepsilon)]$.

Suppose that $d[x(t), \widehat{o}] \geq \delta, \forall t \in [t_0, t_0 + T(\varepsilon)]$.

Define $m(t) = V(t, x(t)) + \int_{t_0}^t V^*(\tau, x(\tau))d\tau, t \geq t_0$.

By the same argument in the proof of Theorem 3.3, we have $m(t) \leq r(t; t_0, w_0), t \geq t_0$, where $w_0 = V(t_0, x_0)$ and $r(t; t_0, w_0)$ be the maximal solution of (3.2). Therefore,

$$0 \leq V(t_0 + T, x(t_0 + T))$$
$$\leq r(t_0 + T; t_0, w_0) - \int_{t_0}^{t_0+T} V^*(\tau, x(\tau))d\tau$$
$$\leq r(t_0 + T; t_0, w_0) - Tc(\delta)$$

Since $V(t_0, x_0) \leq b(d[x_0, \widehat{o}]) < b(\delta_0)$, so we have $w_0 = V(t_0, x_0) < b(\delta_0)$ and hence $r(t_0+T; t_0, w_0) < a(\rho)$. Therefore, $0 \leq V(t_0 + T, x(t_0 + T)) < a(\rho) - Tc(\delta) < 0$. This contradiction shows that there exists $t_1 \in [t_0, t_0 + T]$ such that $d[x(t_1), \widehat{o}] < \delta$. On the other hand, $x(t; t_1, x(t_1; t_0, x_0)) = x(t; t_0, x_0), \forall t \geq t_1$, hence,

$$d[x(t), \widehat{o}] < \varepsilon, \forall t \geq t_0 + T(\varepsilon).$$

This shows that the trivial solution $x = \widehat{o}$ of (3.1) is uniformly asymptotically stable. The proof is completed. □

**Theorem 3.5.** *Suppose that there exists a function $V(t, x)$ satisfies:*

(i) $|V(t, x) - V(t, y)| \leq L(t)d[x, y], \forall (t, x), (t, y) \in \mathbb{R}_+ \times S(\rho)$;

(ii) $\lambda(d[x, \widehat{o}])^p \leq V(t, x) \leq \Lambda(d[x, \widehat{o}])^q$

(iii) $D_f^+ V(t, x) \leq -\gamma(d[x, \widehat{o}])^q + Ke^{-\delta t}, t \geq 0$; *where* $\lambda, \Lambda, \gamma, K, p, q, \delta$ *are positive numbers.*

*If* $\delta > \dfrac{\gamma}{\Lambda} > 0$ *then the trivial solution* $x = \widehat{o}$ *of (3.1) is uniformly exponentialy stable.*

*Proof.* By Theorem 3.2, $x = \widehat{o}$ is uniformly stable. Hence, there exists $H > 0$ such that if $t_0 \in \mathbb{R}_+$ and $d[x_0, \widehat{o}] < H$ then $d[x(t; t_0, x_0), \widehat{o}] < \rho, \forall t \geq t_0$.

Let $M = \dfrac{\gamma}{\Lambda}, m(t) = V(t, x(t))e^{M(t-t_0)}, t \geq t_0$. We have

$$d^+m(t) \leq MV(t, x(t))e^{M(t-t_0)} + e^{M(t-t_0)}D_f^+ V(t, x(t))$$
$$\leq MV(t, x(t))e^{M(t-t_0)} + e^{M(t-t_0)}[Ke^{-\delta t} - \gamma(d[x, \widehat{o}])^q]$$
$$\leq MV(t, x(t))e^{M(t-t_0)} + Ke^{(M-\delta)(t-t_0)} -$$
$$- \dfrac{\gamma}{\Lambda}e^{M(t-t_0)}V(t, x(t)) = Ke^{(M-\delta)(t-t_0)}.$$

Apply Lemma 3.1, $m(t) - m(t_0) \leq K\int_{t_0}^t e^{(M-\delta)(\tau-t_0)}d\tau = \dfrac{K}{M-\delta}[e^{(M-\delta)(t-t_0)} - 1]$. By hypotheses, $m(t_0) = V(t_0, x_0) \leq \Lambda(d[x_0, \widehat{o}])^q$, we have

$$m(t) \leq \dfrac{K}{M-\delta}e^{(M-\delta)(t-t_0)} - \dfrac{K}{M-\delta} + \Lambda(d[x_0, \widehat{o}])^q$$

Put $\delta_1 = -(M - \delta) > 0$, then

$$m(t) \leq \Lambda(d[x_0, \widehat{o}])^q + \dfrac{K}{\delta_1} - \dfrac{K}{\delta_1}e^{-\delta_1(t-t_0)} \leq \Lambda(d[x_0, \widehat{o}])^q + \dfrac{K}{\delta_1}, t \geq t_0.$$

Therefore $V(t, x(t)) \leq \beta_1(d[x_0, \widehat{o}])e^{-M(t-t_0)}, t \geq t_0$, where $\beta_1(d[x_0, \widehat{o}]) = \Lambda(d[x_0, \widehat{o}])^q + \dfrac{K}{\delta_1}$. On the other hand, $\lambda(d[x(t), \widehat{o}])^p \leq V(t, x(t)), t \geq t_0$, we have

$$d[x(t), \widehat{o}] \leq [\dfrac{\beta_1(d[x_0, \widehat{o}])}{\lambda}]^{\frac{1}{p}} e^{-\frac{M}{p}(t-t_0)}, t \geq t_0$$

Denote $\alpha = \dfrac{M}{p}, \beta(d[x_0, \widehat{o}]) = [\dfrac{\beta_1(d[x_0, \widehat{o}])}{\lambda}]^{\frac{1}{p}}$ then

$$d[x(t), \widehat{o}] \leq \beta(d[x_0, \widehat{o}])e^{-\alpha(t-t_0)}, t \geq t_0$$

This shows that the trivial solution $x = \widehat{o}$ of (3.1) is uniformly exponentially stable. The proof is completed. □

MATHEMATICS DEPARTMENT, HA NOI UNIVERSITY OF EDUCATION
136 XUAN THUY ROAD, CAU GIAY DISTRICT, HA NOI - VIET NAM
*E-mail address: Hienlv@dhsphn.edu.vn*